\documentclass[
 aip,
 amsmath,amssymb,
 reprint,
]{revtex4-1}

\usepackage{graphicx}
\usepackage{bm}
\usepackage{hyperref}
\usepackage{subcaption}
\usepackage{mathtools}

\begin{document}

\title{A conservation-consistent boundary condition for nonlinear models of soluble-surfactant-laden falling films}

\author{Sanghasri Mukhopadhyay}
\affiliation{Centre for Applied Sciences (CAS), International Institute of Information Technology, Bengaluru 560100, India}

\author{Severine Millet}
\affiliation{Universite Claude Bernard Lyon 1, LMFA, UMR5509, CNRS, Ecole Centrale de Lyon, INSA Lyon, 69622 Villeurbanne France}

\author{Bastien Di Pierro}
\affiliation{Universite Claude Bernard Lyon 1, LMFA, UMR5509, CNRS, Ecole Centrale de Lyon, INSA Lyon, 69622 Villeurbanne France}

\author{Asim Mukhopadhyay}
\affiliation{Department of Mathematics, Vivekananda Mahavidyalaya,
Burdwan 713103, West Bengal, India}

\begin{abstract}

A conservation-consistent boundary condition is proposed for nonlinear models of soluble-surfactant-laden falling films, ensuring exact conservation of total surfactant mass. The formulation resolves an inconsistency in widely used reduced models \cite{pascalStabilityInclinedFlow2019, dalessioMarangoniInstabilitiesAssociated2020}, which exhibit a gradual drift of mass during nonlinear evolution in a closed periodic domain. We show that this originates from an inconsistency in the surface transport reduction and derive a corrected boundary condition that removes this defect. As the discrepancy appears only at nonlinear order, linear stability results remain unaffected, explaining why the issue has remained unnoticed\cite{samantaEffectSolubleSurfactant2025, samantaRoleSolubleSurfactant2025}.

\end{abstract}

\maketitle



Gravity-driven thin liquid films containing soluble surfactants exhibit a wide range of interfacial instabilities arising from the interplay between hydrodynamics, surface tension, and Marangoni stresses. Reduced-order models based on long-wave and weighted-residual approaches have therefore become essential tools for studying the nonlinear evolution of such systems. Among these, the formulations developed by Pascal et al. (2019)\ \cite{pascalStabilityInclinedFlow2019} and later by D'Alessio et al. (2020)\ \cite{dalessioMarangoniInstabilitiesAssociated2020} are widely used for soluble-surfactant-laden falling films.

While reproducing the nonlinear simulations reported in these works, we observed a slow but systematic drift in the total surfactant mass during time evolution. In a closed system without external sources or sinks, such behaviour is unphysical and indicates an inconsistency in the reduced model. A closer examination reveals that the issue originates from the boundary condition governing surfactant exchange between the bulk and the free surface. Although Pascal et al. (2019)\ \cite{pascalStabilityInclinedFlow2019} noted the lack of strict mass conservation, the underlying cause remained unclear. Here, we revisit the surface transport reduction and derive a conservation-consistent boundary condition that restores exact global mass conservation.


Pereira and Kalliadasis (2008)\cite{pereiraTransportEquationInterfacial2008} derived a general transport equation for interfacial quantities such as surfactants and applied it to insoluble surfactant-laden falling films \cite{pereiraDynamicsFallingFilm2008}. Pascal et al. (2019)\ \cite{pascalStabilityInclinedFlow2019} and D’Alessio et al. (2018)\ \cite{dalessioMarangoniInstabilitiesAssociated2020} extended this framework to the soluble case by incorporating an adsorption-desorption source term \(J_{bs}\). The three-dimensional form of the above transport equation with $J_{bs}$ terms is as follows:
\begin{equation}
\partial_t \Gamma
+\vec{V}\cdot \boldsymbol{\nabla}_{xy}\Gamma
+\Gamma\,\boldsymbol{\nabla}_s\cdot \vec{V}
=
D_s \nabla_s^2 \Gamma + J_{bs}.
\label{wrong_bc}
\end{equation}
In the above equation, $x-$axis represents the streamwise direction, the $y-$axis represents the spanwise direction and the $z-$axis  represents the cross-stream direction.  $\Gamma(x,y,t)$ is the surface surfactant concentration, $\vec{V}$ is the surface velocity, $\boldsymbol{\nabla}_{xy}$ is the in-plane gradient operator, $\boldsymbol{\nabla}_s=(I-\mathbf{nn})\cdot\boldsymbol{\nabla}$ is the surface gradient operator, $D_s$ is the surface diffusivity, and $J_{bs}$ denotes the adsorption-desorption flux.
 \\
In the present work, we instead adopt the conservative surface-balance form used in earlier studies of interfacial surfactant transport (Gaver and Grotberg (1990)\cite{gaverDynamicsLocalizedSurfactant1990}, Stone (1990) \cite{stoneSimpleDerivationTimedependent1990}, Manikantan and Squires (2020)\cite{manikantanSurfactantDynamicsHidden2020}), namely
\begin{equation}
\underbrace{\partial_t \Gamma}_\text{unsteady} + \underbrace{\boldsymbol{\nabla}_s\cdot(\Gamma \vec{V}_s)}_\text{advection}
=
\underbrace{D_s \nabla_s^2 \Gamma}_\text{diffusion} + \underbrace{J_{bs},}_\text{source term}
\label{right_bc}
\end{equation}
where $\vec{V}_s=(I-\mathbf{nn})\cdot \vec{V}$ is the tangential velocity on the free surface. This form preserves the conservative structure of the surface transport equation and, when reduced consistently, restores the global surfactant mass balance. 
We want to emphasize that $\nabla_s \cdot(\Gamma \vec{V}_s)$ is the surface advection term, while $D_s \nabla_s^2 \Gamma$ is the surface diffusion term. \par
It should be noted that in writing equation (\ref{right_bc}), we have neglected the contribution of the surface divergence ${\boldsymbol{\nabla}_s\cdot(\Gamma \vec{V})}$ \cite{stoneSimpleDerivationTimedependent1990} for the normal component to the surface advective flux $\mathbf{n}(\mathbf{n}\cdot \vec{V}) \Gamma$ \cite{manikantanSurfactantDynamicsHidden2020} (where $\vec{V}$ is the full velocity of the free surface and $\mathbf{n}$ the unit normal), which is $\Gamma ({\boldsymbol{\nabla}_s}\cdot \mathbf{n})(\vec{V} \cdot \mathbf{n})$.  This is because the term $\Gamma ({\boldsymbol{\nabla}_s}\cdot \mathbf{n})(\vec{V} \cdot \mathbf{n})$ is a source-like contribution that accounts for variations in surface concentration and is simply related to the product of mean curvature and normal velocity \cite{stoneSimpleDerivationTimedependent1990}, and mean curvature is extremely small in order of magnitude ($O(\varepsilon^2))$. In addition, the Langmuir term $J_{bs}$ is used as a source term for soluble surfactants. We would like to highlight that in this model the surfactant is permitted to diffuse along the interface with constant diffusivity or to be advected along the thin film's surface \cite{gaverDynamicsLocalizedSurfactant1990, stoneSimpleDerivationTimedependent1990, rahmanNanoscaleSurfactantTransport2025}.

For a two-dimensional flow ($v=0$, $\partial_y=0$), Eq.~\eqref{right_bc} reduces in dimensionless form to
\begin{eqnarray}
\varepsilon \left[
\frac{\partial \Gamma}{\partial t}
+\frac{1}{\sqrt{1+\varepsilon^2 h_x^2}}
\frac{\partial(\Gamma u_s)}{\partial x}
\right]
\nonumber\\
=
\frac{\varepsilon^2}{Pe_s\sqrt{1+\varepsilon^2 h_x^2}}
\frac{\partial}{\partial x}
\left[
\frac{1}{\sqrt{1+\varepsilon^2 h_x^2}}
\frac{\partial \Gamma}{\partial x}
\right]
\nonumber\\
+k_s[\kappa(1-\Gamma)C-\Gamma],
\label{write_bc_surface_surfactant}
\end{eqnarray}
where, $\varepsilon$ is the aspect ratio, $C(x,z,t)$ and $\Gamma(x,t)$ denote the bulk and interfacial surfactant concentrations, respectively, $Pe_s$ is the surface P\'eclet number, $k_s$ is the desorption rate, $\kappa$ is the adsorption equilibrium constant and

\[
u_s(x,t)=\frac{u(x,z=h,t)}{\sqrt{1+\varepsilon^2 h_x^2}}
\]
is the surface velocity. All the detail definitions of different parameters can bo found in the supplementary file. Equation \eqref{write_bc_surface_surfactant} differs from the corresponding closure used in Pascal et al. (2019)\cite{pascalStabilityInclinedFlow2019} and D'Alessio et al. (2020)\cite{dalessioMarangoniInstabilitiesAssociated2020} (equation (2.6) in their paper). They have an additional
term, $\Gamma \partial_x h\partial_z u$ up to first order in $\varepsilon$, which is non-zero due to the correction for soluble surfactant in the assumed specific velocity profile that is necessary for the weighted residual method. This term is critical in understanding the increase in mass on the free surface.
In contrast, equation \eqref{write_bc_surface_surfactant}  closely resembles Pereira and Kalliadasis’ conservative surface-transport form\cite{pereiraDynamicsFallingFilm2008} in the insoluble-surfactant limit.

The detailed derivation of Eq.~\eqref{write_bc_surface_surfactant} and its reduction in curvilinear coordinates are given in the Supplementary Material.\\[10pt]
\begin{figure}
\centering
\includegraphics[width=0.8\linewidth]{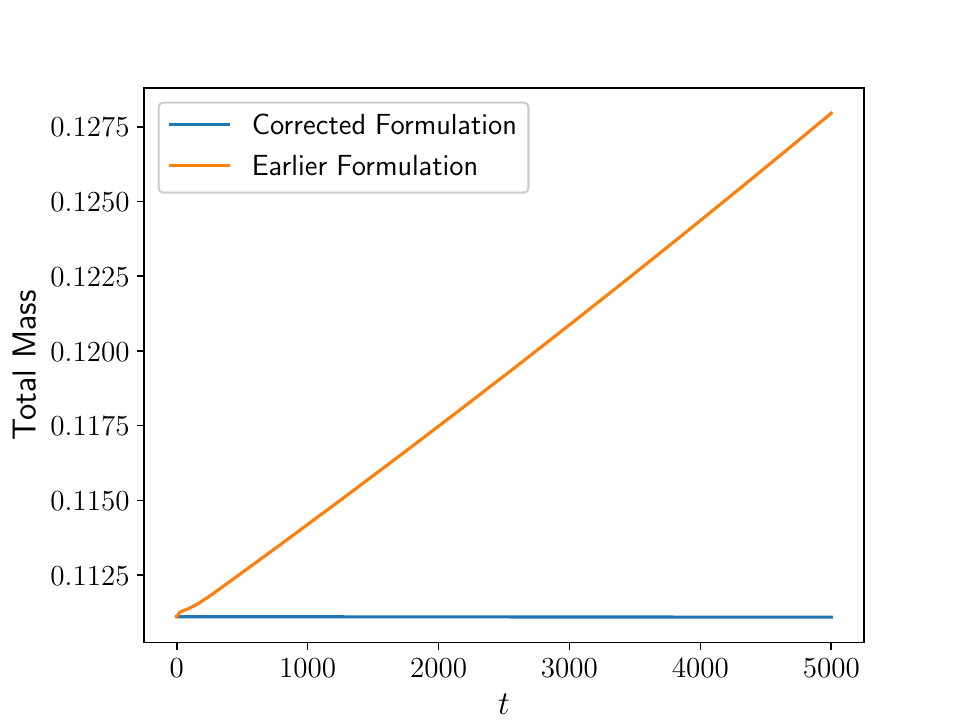}
\caption{
Temporal evolution of the total surfactant mass $M(t)$ for the earlier formulation of  Pascal et al. (2019)\ \cite{pascalStabilityInclinedFlow2019} and D'Alessio et al. (2020)\ \cite{dalessioMarangoniInstabilitiesAssociated2020}  and the new conservation-consistent model.}
\label{fig_pascal}
\end{figure}
\begin{figure*}
    \centering
    \subcaptionbox{earlier formulation}{\includegraphics[width=0.49\textwidth]{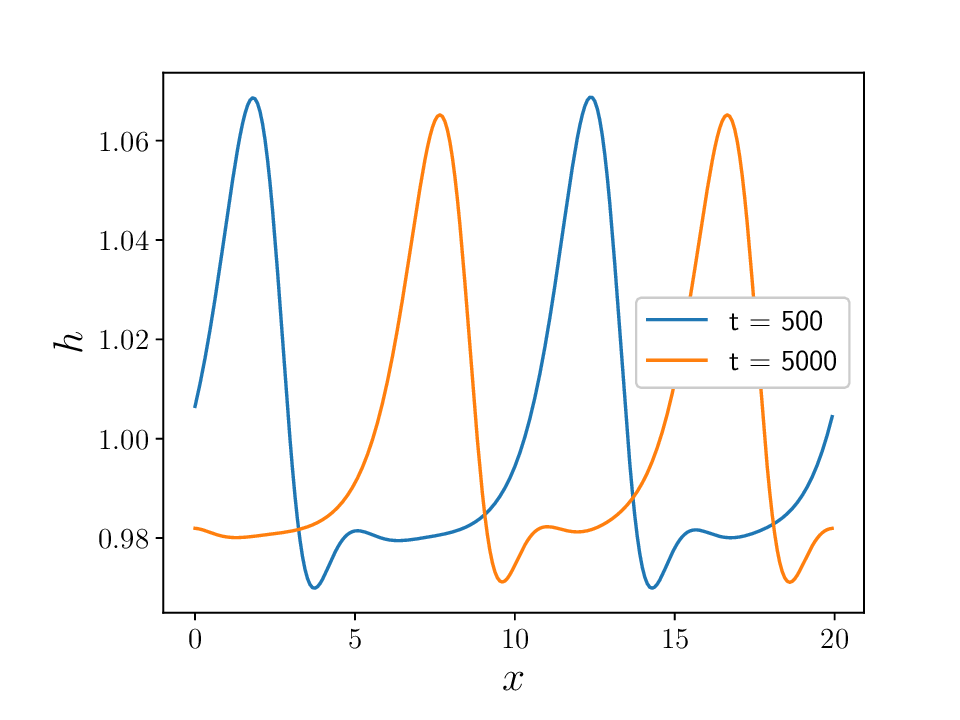}}
    \hfill
    \subcaptionbox{corrected formulation}{\includegraphics[width=0.49\textwidth]{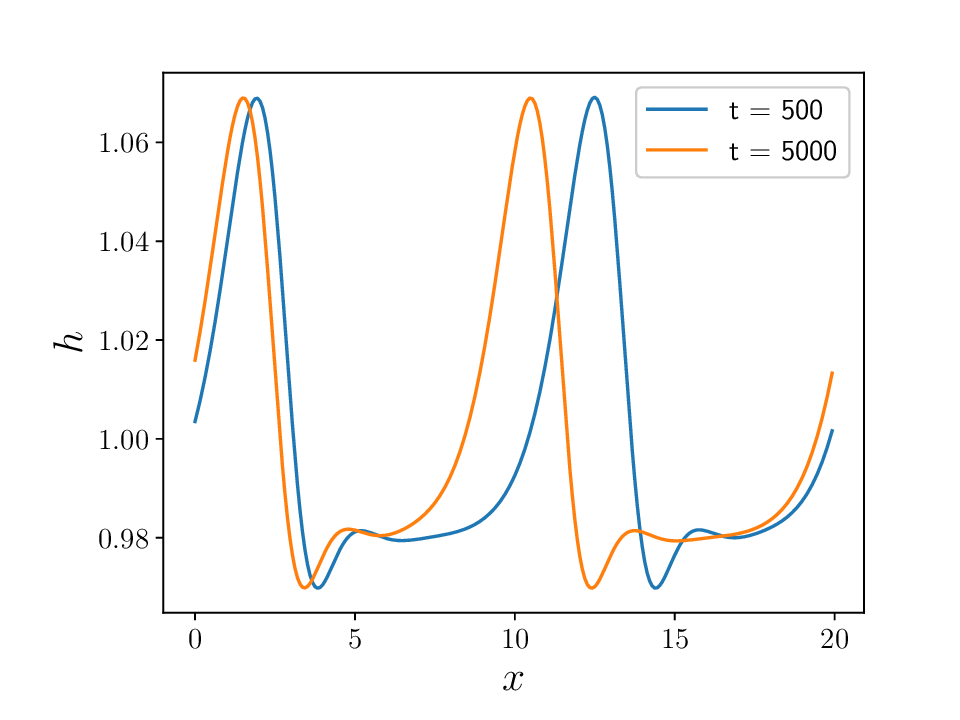}}

    \subcaptionbox{earlier formulation}{\includegraphics[width=0.49\textwidth]{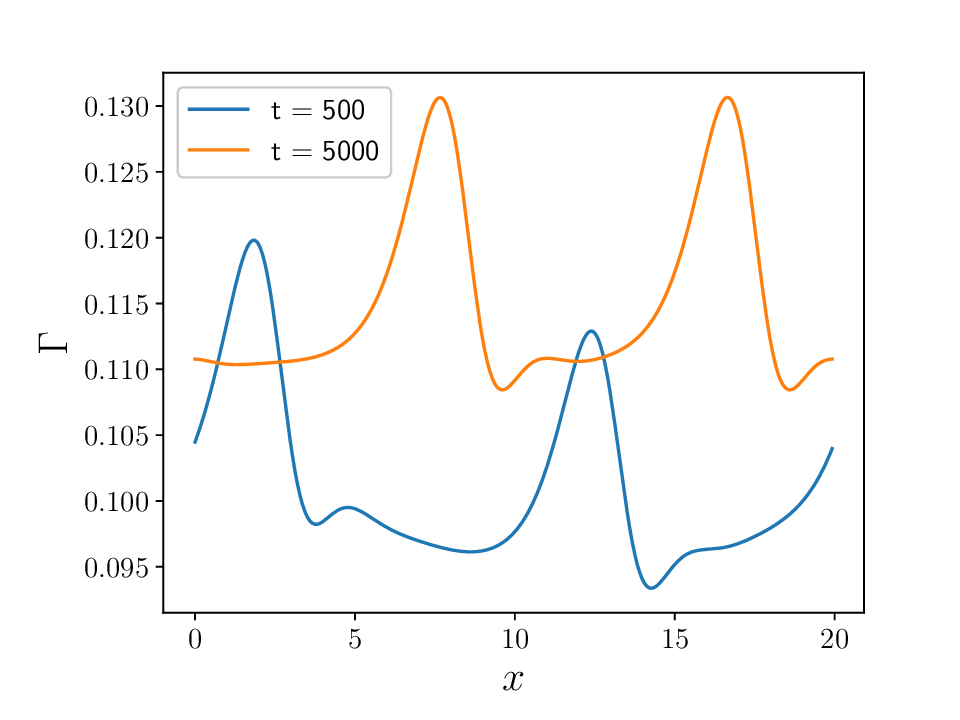}}
    \hfill
    \subcaptionbox{corrected formulation}{\includegraphics[width=0.49\textwidth]{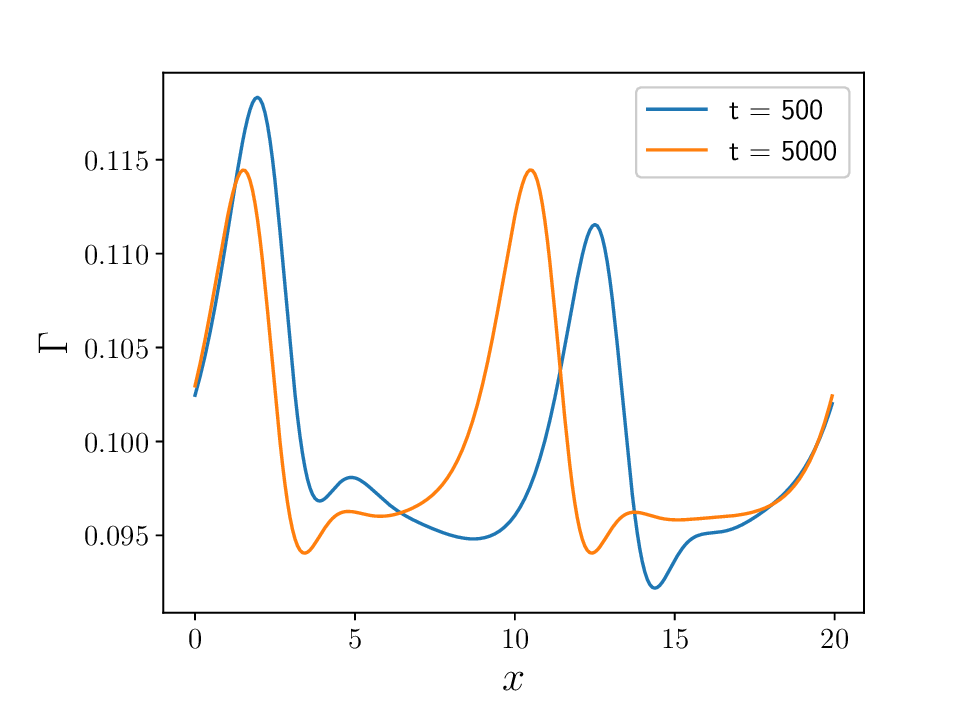}}

    \subcaptionbox{earlier formulation}{\includegraphics[width=0.49\textwidth]{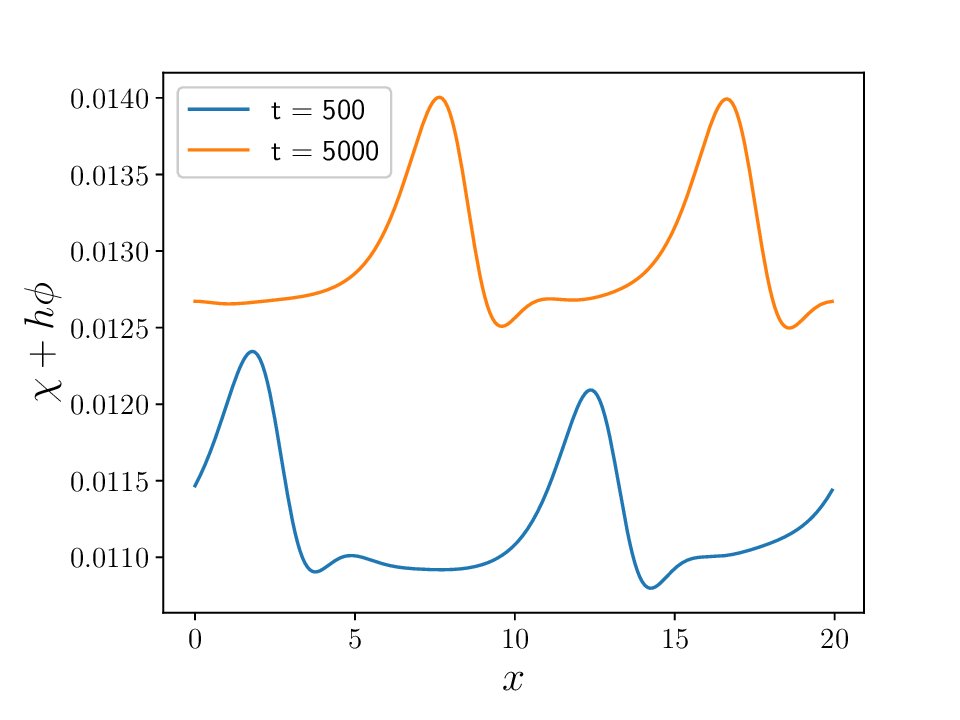}}
    \hfill
    \subcaptionbox{corrected formulation}{\includegraphics[width=0.49\textwidth]{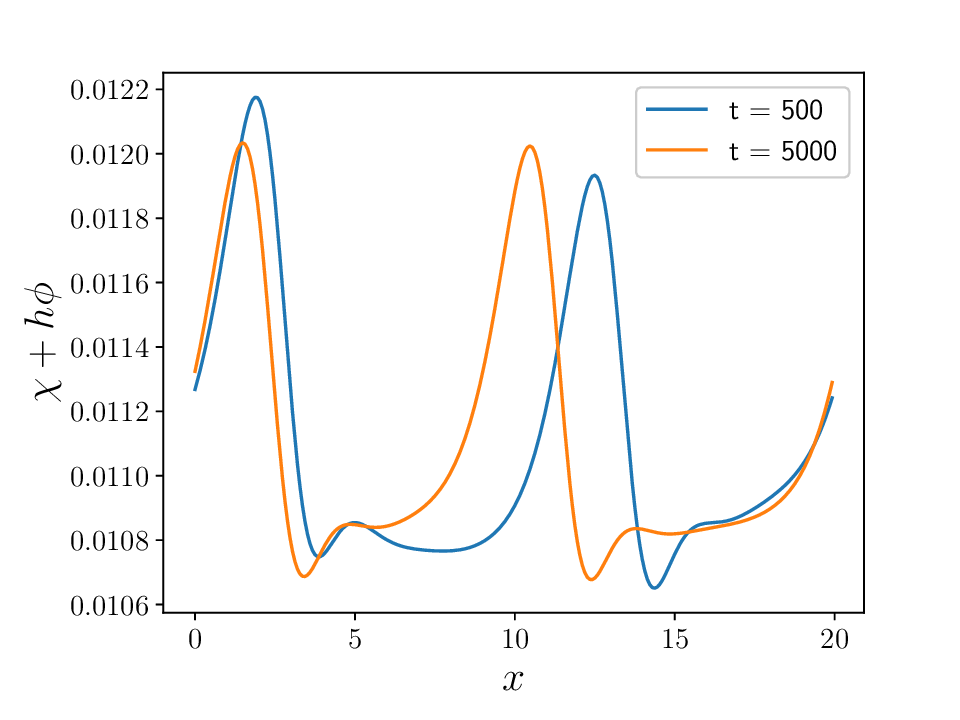}}

\caption{
Comparison of nonlinear evolution for the earlier formulation (left column) and the corrected formulation (right column). Shown are the film thickness $h(x,t)$ (top row), surface concentration $\Gamma(x,t)$ (middle row), and bulk contribution $\chi + h\phi$ (bottom row) at $t=500$ and $t=5000$.
} 
    \label{fig:model_comp}
\end{figure*}
\begin{figure*}
    \centering
    \subcaptionbox{earlier formulation}{\includegraphics[width=0.49\textwidth]{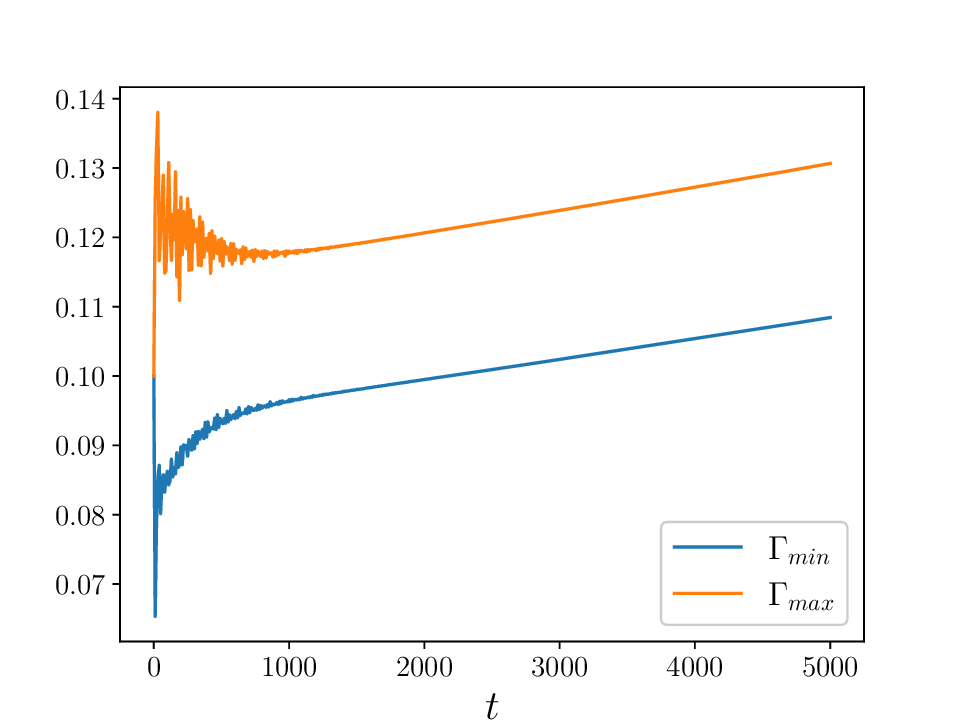}}
    \hfill
    \subcaptionbox{corrected formulation}{\includegraphics[width=0.49\textwidth]{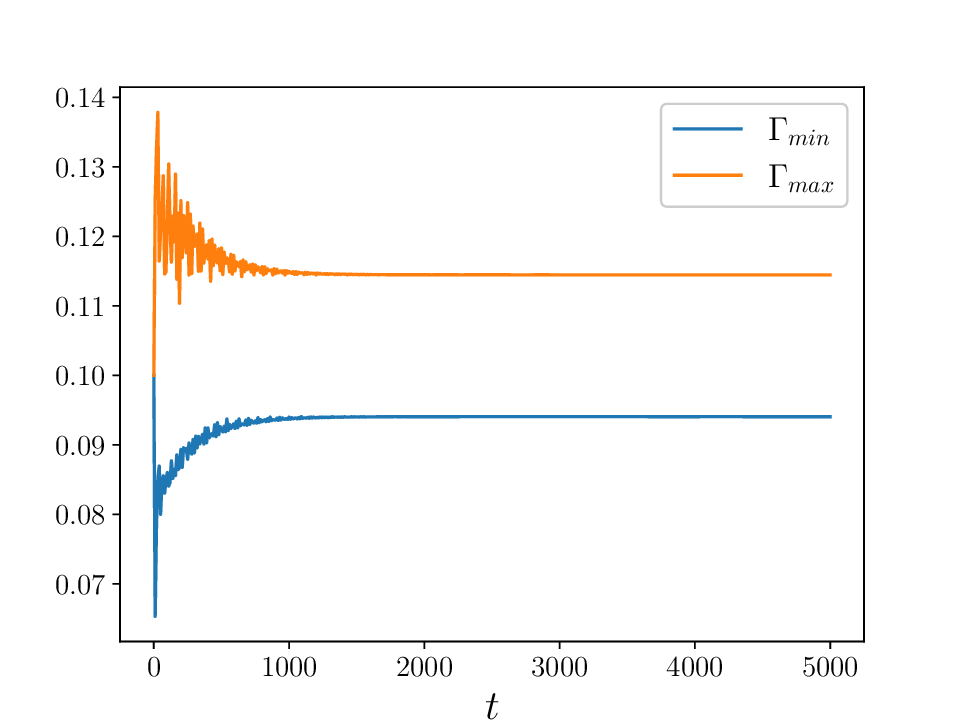}}
   \caption{
Temporal evolution of the minimum and maximum surface concentration, $\Gamma_{\min}$ and $\Gamma_{\max}$, for the earlier formulation (left) and the conservation-consistent model (right).}
    \label{fig:gamma_comp}
\end{figure*}

To assess its impact, we reproduced the nonlinear simulations of Pascal et al. (2019)\ \cite{pascalStabilityInclinedFlow2019} and D'Alessio et al. (2020)\ \cite{dalessioMarangoniInstabilitiesAssociated2020} and monitored the total surfactant mass
\begin{equation}
M(t)=\int_0^L \int_0^{h(x,t)} C(x,z,t)\,dz\,dx+\int_0^L \Gamma(x,t)\,dx,
\label{total_mass}
\end{equation}
where, $L$ is the arbitrary domain length being considered in a stream-wise direction. For a closed periodic system without external sources or sinks, $M(t)$ must remain constant in time. 

In terms of the reduced variables used in the weighted-residual model, this can be written as
\begin{equation}
M(t)=\int_0^L \left[\phi(x,t)h(x,t)+\chi(x,t)+\Gamma(x,t)\right]dx,
\label{total_mass_reduced}
\end{equation}
where $\phi(x,t)$ represents the leading-order bulk concentration evaluated at the interface and $\chi(x,t)$ is the excess/deficit amount of surfactant in the bulk. The detailed expressions for $\phi$ and $\chi$, as well as the full mass-balance diagnostics, are provided in the Supplementary Material.

Figure~\ref{fig_pascal} shows the evolution of the total surfactant mass obtained with the earlier formulation and with the present one. All simulations are performed for $Re = 1.5$, $Fr = 0.7071$, $Pe_b = Pe_s = 700$, $\varepsilon = 0.1$, $Mr = 1$, $k_s = 1$, $\kappa = 10$, $\Gamma_e = 0.1$, and $Ka = 0.75$ in a periodic domain of length $L = 20$. The earlier model (old model) exhibits a systematic drift, whereas the corrected formulation (new model) preserves the total mass to numerical accuracy.

The nonlinear dynamics themselves are also modified. Figure~\ref{fig:model_comp} compares the wave profile, surface concentration, and bulk surfactant contribution obtained using the earlier coefficient and the corrected one. The difference can be traced to the nonlinear coupling term in the reduced surface transport equation. In the earlier formulation, the term
$$\Gamma\,\partial_x\Gamma\,\partial_x h$$
appears with coefficient $5/4$, whereas the conservation-consistent reduction yields the coefficient $1/4$. Although small in appearance, this difference is decisive in the nonlinear dynamics.

For completeness, the full second-order weighted-residual model is given below.
\begin{widetext}
\begin{equation}\label{mass_conservation}
    \frac{\partial h}{\partial t} + \frac{\partial q}{\partial x} = 0
\end{equation}

 \begin{eqnarray}\label{model_q}
    \frac{\partial q}{\partial t} + \frac{\partial}{\partial x} \left[ \frac{9}{7} \frac{q^2}{h} + \frac{5}{12} \frac{Re}{Fr^2}\frac{\cot \beta}{Re}h^2 + \frac{5}{4} Mr \Gamma\right]= \frac{1}{7}\frac{q}{h}\frac{\partial q}{\partial x}  +\frac{5}{6}\frac{Ka}{Re}h \frac{\partial^3 h}{\partial x^3}
    \nonumber\\
    + \frac{1}{\varepsilon Re} \frac{5}{2} \left(\frac{1}{3}\frac{Re}{Fr^2}h - \frac{q}{h^2} \right)  + \frac{\varepsilon}{Re} \left[\frac{9}{2} \frac{\partial^2 q}{\partial x^2}  - \frac{9}{2} \frac{1}{h} \frac{\partial q}{\partial x} \frac{\partial h}{\partial x} + 4 \frac{q}{h^2} \left(\frac{\partial h}{\partial x} \right)^2 -6 \frac{q}{h} \frac{\partial^2 h}{\partial x^2}\right] 
    \nonumber\\
    +\frac{\varepsilon Re Mr}{16} \left[ \frac{1}{3}h^2 \frac{\partial^2 \Gamma}{\partial t \partial x} + \frac{15}{14}hq \frac{\partial^2 \Gamma}{\partial x^2} +\frac{19}{21}h \frac{\partial q}{\partial x} \frac{\partial \Gamma}{\partial x} + \frac{5}{7}q \frac{\partial \Gamma}{\partial x}\frac{\partial h}{\partial x}\right] 
\end{eqnarray}
\begin{eqnarray}\label{model_chi}
    \frac{\partial}{\partial t}
\left( \chi + \phi h \right)
+ \frac{\partial}{\partial x} \left[
  \frac{33}{40} \frac{q\,\chi}{h}
+ \phi q
\right] = \frac{\varepsilon}{\mathrm{Pe}_b} \left[
\frac{\partial^2 \chi}{\partial x^2}
- \frac{3\chi}{h^2} \left( \frac{\partial h}{\partial x} \right)^2
+ h \frac{\partial^2 \phi}{\partial x^2}
\right] 
\nonumber\\ 
\underbrace{- \frac{3}{\varepsilon Pe_b} \frac{\chi}{h^2}}_{\text{source term}}  
- \frac{\mathrm{3}}{80} \varepsilon Re Mr
\left[
\frac{\partial^2 \Gamma}{\partial x^2} h \chi
+ \frac{\partial \chi}{\partial x} \frac{\partial \Gamma}{\partial x} h
+ \chi \frac{\partial \Gamma}{\partial x} \frac{\partial h}{\partial x}
\right]
\end{eqnarray}
\begin{equation}\label{model_Gamma}
\begin{aligned}
\frac{\partial \Gamma}{\partial t} 
+ \frac{3}{2}\frac{\partial}{\partial x} \left(\frac{q \Gamma}{h}\right)
&= \frac{\varepsilon}{Pe_s} \frac{\partial^2 \Gamma}{\partial x^2} 
\underbrace{+ \frac{3}{\varepsilon Pe_b}\frac{\chi}{h^2}}_{\text{adsorption term}} 
+ \frac{1}{4} \varepsilon Re Mr \Bigg[
\frac{\partial^2 \Gamma}{\partial x^2} \Gamma h 
+ \left(\frac{\partial \Gamma}{\partial x}\right)^2 h 
+ \underbrace{\Gamma\frac{\partial \Gamma}{\partial x}\frac{\partial h}{\partial x}}_{\text{coupling term}}
\Bigg].
\end{aligned}
\end{equation}
\end{widetext}
The derivation of the reduced model and the momentum equation are given in the Supplementary Material. For equilibrium of the surfactant concentration, the coefficient of $\chi/h^2$ in the bulk and surface transport equations must be equal and opposite. In the formulations of Pascal et al. (2019)\cite{pascalStabilityInclinedFlow2019} and D'Alessio et al. (2020)\cite{dalessioMarangoniInstabilitiesAssociated2020}  this condition is not satisfied, owing to a different definition of $\chi$ involving a distinct numerical coefficient.

Figure~\ref{fig:gamma_comp} compares the extrema of the surface concentration in time for the earlier and corrected formulations. The earlier model shows secular, unphysical growth, while the corrected model remains bounded and conservation-consistent.

Although the difference between the two boundary conditions appears minor, it leads to a fundamental violation of global surfactant conservation in nonlinear simulations. In particular, the coefficient change from $5/4$ to $1/4$ in the coupling term $\Gamma\,\partial_x \Gamma\,\partial_x h$ induces unphysical growth of the surfactant mass in the earlier formulation. Since this discrepancy arises only at nonlinear order, it does not affect the linearized problem, which explains why it remained unnoticed in previous studies focused on linear stability. The corrected boundary condition therefore establishes a physically consistent and reliable framework for nonlinear investigations of soluble-surfactant falling films.

\bibliographystyle{aipnum4-1}
\bibliography{bibliography}

\begin{widetext}
\section*{Supplementary Material}

The Supplementary material contains the detailed derivation of the conservative surface transport equation, the full formulation of the weighted-residual model, and additional diagnostics of bulk, interfacial, and total surfactant mass fluxes supporting the results presented in the main text. 
\end{widetext}

\newpage \null \newpage

\section*{S1. Governing equations and boundary conditions}

Consider a two-dimensional flow of a viscous incompressible fluid film with density $\rho$, viscosity $\mu$ in an inclined plane of inclination $\theta$ with the horizon under the action of gravity, the $x$ and $z$ axes are chosen along and normal to the flow direction, respectively,  $u(t,x,z)$ and $w(t,x,z)$ denote the velocity components in the $x$ and $z$ directions, respectively, $t$ denotes the time, $p(t,x,z)$ denotes the pressure, $\textsl{g} = (g \sin\theta, -g\cos\theta)$ is the gravitational acceleration, and $C(t,x,z)$ is the bulk concentration of the surfactant. The free surface (fluid-air interface) is given by $z = h(x,t)$. The liquid contains a fixed amount of soluble surfactant. Since we are working with a soluble surfactant, we must consider two different concentrations: one for the
surfactant adsorbed at the surface and one for the surfactant dissolved in the bulk of the layer. Furthermore, we will assume that the concentrations of surfactants are below the critical level for micelle formation, such that the surfactant added to the liquid exists as monomers. The governing equations in dimensionless form are
\begin{widetext}
\begin{equation}\label{mass}
	\dfrac{\partial u}{\partial x} + \dfrac{\partial w}{\partial z} = 0,
\end{equation}
\begin{eqnarray}\label{x-momentum}
	\varepsilon Re\left(\dfrac{\partial u}{\partial t} + u \dfrac{\partial u}{\partial x} + w \dfrac{\partial u}{\partial z}\right) = - \varepsilon Re \dfrac{\partial p}{\partial x}+ \frac{Re}{Fr^2} + \left(\varepsilon^2 \dfrac{\partial^2 u}{\partial x^2} + \dfrac{\partial^2 u}{\partial z^2}\right),
\end{eqnarray}
\begin{eqnarray}\label{y-momentum}
	\varepsilon^2 Re\left(\dfrac{\partial w}{\partial t} + u \dfrac{\partial w}{\partial x} + w \dfrac{\partial w}{\partial z}\right)  = -  Re \dfrac{\partial p}{\partial z} - \frac{Re}{Fr^2} \cot \theta  
	+ \left(\varepsilon^3 \dfrac{\partial^2 w}{\partial x^2} + \varepsilon \dfrac{\partial^2 w}{\partial z^2}\right),
\end{eqnarray}
\begin{equation}\label{concentration}
	\varepsilon Pe_b\left(\dfrac{\partial C}{\partial t} + u \dfrac{\partial C}{\partial x} + w \dfrac{\partial C}{\partial z}\right)  = \varepsilon^2 \dfrac{\partial^2 C}{\partial x^2} + \dfrac{\partial^2 C}{\partial z^2}
\end{equation}
\text{The boundary conditions at} $z = 0$ \text{are as follows:}
\begin{equation}\label{boundary_bottom}
	u = 0, \qquad w = 0, \qquad  \dfrac{\partial C}{\partial z} = 0,
\end{equation}
\text{The boundary conditions at} $z = h$ \text{are given by:} 
\begin{eqnarray}\label{boundary_normal}
	p_a - p + 2 Re^{- 1}\varepsilon \left\{\varepsilon^2 \dfrac{\partial u}{\partial x} \left(\dfrac{\partial h}{\partial x}\right)^2 - \left(\dfrac{\partial u}{\partial z} + \varepsilon^2  \dfrac{\partial w}{\partial x}\right) \dfrac{\partial h}{\partial x} + \dfrac{\partial w}{\partial z} \right\}
	\left(1 +\varepsilon^2 \left(\dfrac{\partial h}{\partial x}\right)^2\right)^{- 1} = \varepsilon^2  \dfrac{\partial^2 h}{\partial x^2}\frac{We - Mr (\Gamma - \Gamma_e)}{\left\{1 +\varepsilon^2 \left(\dfrac{\partial h}{\partial x}\right)^2 \right\}^{3/2}}, \qquad
\end{eqnarray}
\begin{eqnarray}\label{boundary_tangential}
	\left(\dfrac{\partial u}{\partial z} + \varepsilon^2 \dfrac{\partial w}{\partial x}\right)\left\{1 - \varepsilon^2 \left(\dfrac{\partial h}{\partial x}\right)^2\right\}  - 4 \varepsilon^2 \dfrac{\partial u}{\partial x} \dfrac{\partial h}{\partial x} + \varepsilon Re Mr \dfrac{\partial \Gamma}{\partial x} \left\{1 + \varepsilon^2 \left(\dfrac{\partial h}{\partial x}\right)^2\right\}^{1/2}= 0, \qquad
\end{eqnarray}
\begin{eqnarray}\label{boundary_kinematic}
	w = \dfrac{\partial h}{\partial t} +u \dfrac{\partial h}{\partial x},
\end{eqnarray}

\begin{eqnarray}\label{boundary_surface_sarfactant}
\varepsilon \left[ \frac{\partial \Gamma}{\partial t} + \frac{1}{\sqrt{1 + \varepsilon ^2h_x^2}}  \frac{\partial (\Gamma u_s)}{\partial x} \right] =  \frac{\varepsilon ^2}{Pe_s\sqrt{1 + \varepsilon ^2h_x^2}} \frac{\partial}{\partial x} \left[ \frac{1}{\sqrt{1 + \varepsilon ^2 h_x^2}} \frac{\partial \Gamma}{\partial x} \right] + k_s[\kappa (1 - \Gamma) C - \Gamma]  \qquad
\end{eqnarray}

\end{widetext}

The following transformations are used to transform dimensionless. The upper-scripts star denotes the dimensional variables. We chose a pair of separate characteristic length scales, $\hat{\lambda}$ in the longitudinal direction, whose order is equal to the disturbance wavelength and $\hat{h}$, the mean thickness of the film in the transverse direction. $\hat{V}$ is the velocity scale. 
\begin{widetext}
\begin{eqnarray}\label{dimension}
	x = (2 \pi/\hat{\lambda})x^*, \quad z = (1/\hat{h})z^*, \quad h = (1/\hat{h})h^*, \quad
    u = (1/\hat{V})u^*, \quad w = (\hat{\lambda}/ 2 \pi \hat{h}\hat{V})w^*, \nonumber \\ t = (2 \pi \hat{V}/\hat{\lambda})t^*, \quad p = (1/\rho \hat{V}^2)p^*, \quad C = (\hat{h}/\Gamma_\infty)C^*,\quad \Gamma = (1/\Gamma_\infty)\Gamma^*  
\end{eqnarray} 
\end{widetext}

The nondimensional groups are
\begin{widetext}
\begin{eqnarray*}
Re=\frac{\hat{V}\hat{h}}{\nu},\quad
Fr=\frac{\hat{V}}{\sqrt{g\sin\theta\,\hat{h}}},\quad
Pe_b=\frac{\hat{V}\hat{h}}{D_b},\quad Pe_s=\frac{\hat{V}\hat{h}}{D_s},\quad
We=\frac{\sigma_0}{\rho \hat{V}^2\hat{h}},
\nonumber\\
\quad
Mr=\frac{\sigma_\Gamma\Gamma_\infty}{\rho \hat{V}^2\hat{h}},\quad 
k_s=\frac{k_2\hat{h}}{\hat{V}},\quad
\kappa=\frac{k_1}{k_2\hat{h}},\quad 
\varepsilon=\frac{2\pi\hat{h}}{\hat{\lambda}}\ll 1,\quad  Ka = \varepsilon^2 Re We
\end{eqnarray*} 
\end{widetext}
$\Gamma_\infty$ denotes the maximum packing surface concentration, $D_b$ is the diffusivity of the surfactant in the bulk fluid, $D_s$ is the molecular diffusivity of the surface. We have assumed that the surface tension coefficient will vary linearly with the surface surfactant concentration, which is of the form
$\sigma(\Gamma) = \sigma_0 - \sigma_\Gamma (\Gamma - \Gamma_E),$
where, $\Gamma_E$ is the equilibrium surface surfactant concentration and $\sigma_0 = \sigma(\Gamma_E)$ and $\sigma_\Gamma = - \left.\dfrac{\partial \sigma}{\partial \Gamma}\right|_{\Gamma = \Gamma_E}$. 
\begin{figure*}
    \centering
    \subcaptionbox{Bulk flux, earlier model}{\includegraphics[width=0.49\textwidth]{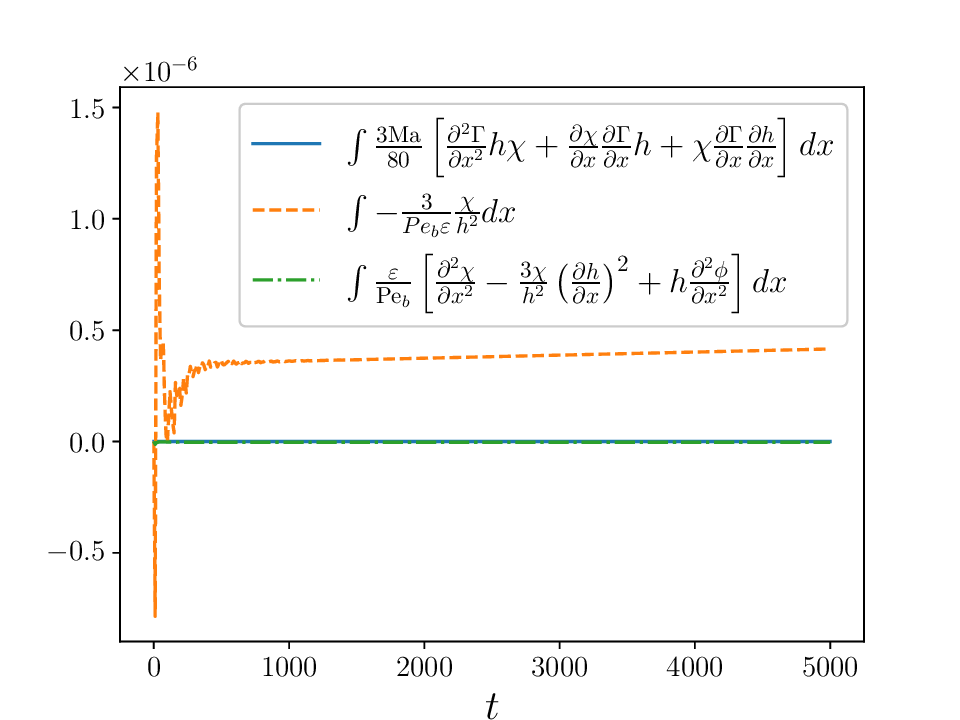}}
    \hfill
    \subcaptionbox{Bulk flux, corrected model}{\includegraphics[width=0.49\textwidth]{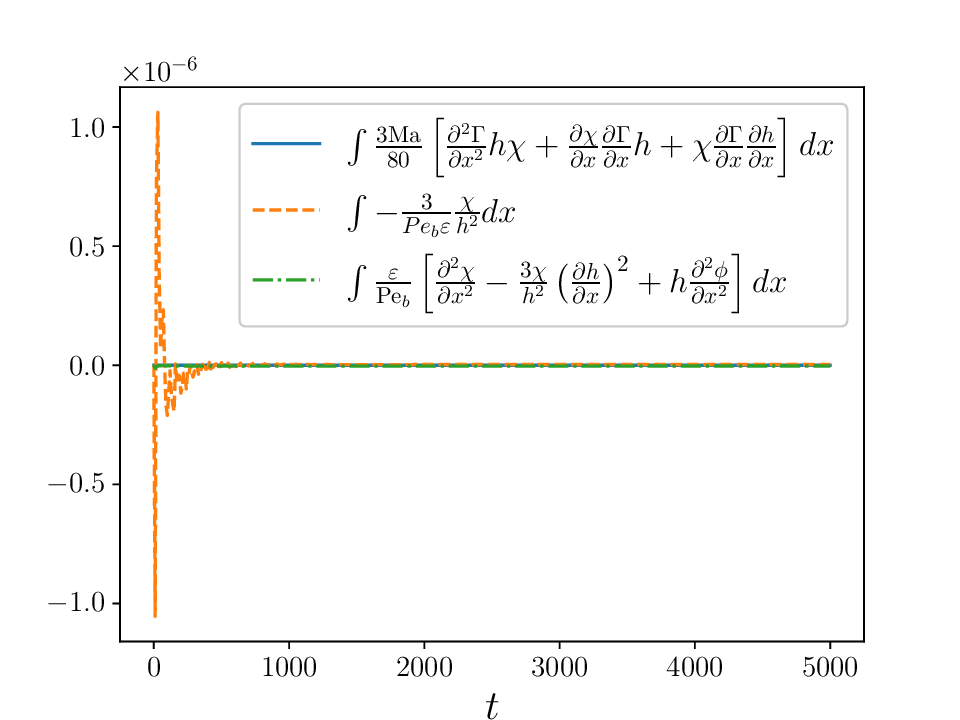}}

    \subcaptionbox{Interfacial flux, earlier model}{\includegraphics[width=0.49\textwidth]{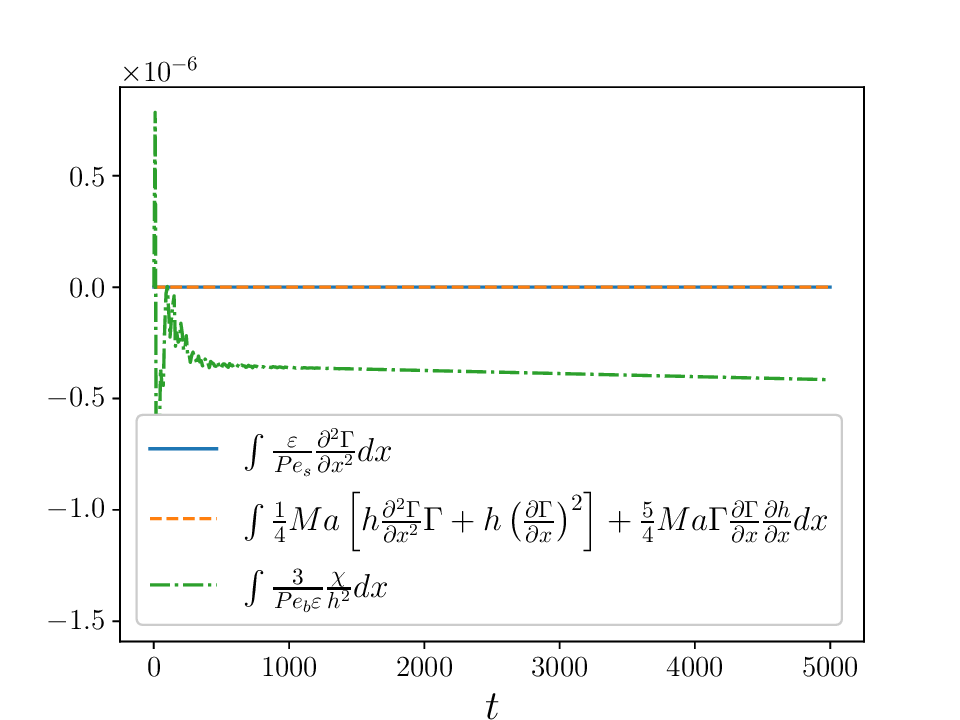}}
    \hfill
    \subcaptionbox{Interfacial flux, corrected model}{\includegraphics[width=0.49\textwidth]{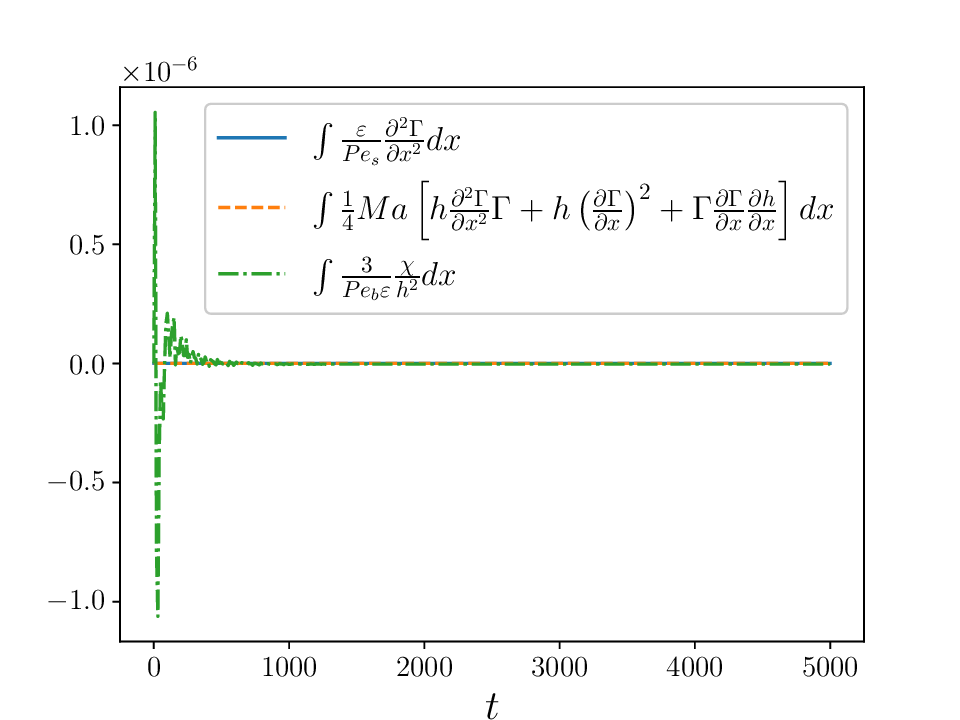}}

    \subcaptionbox{Mass evolution, earlier model}{\includegraphics[width=0.49\textwidth]{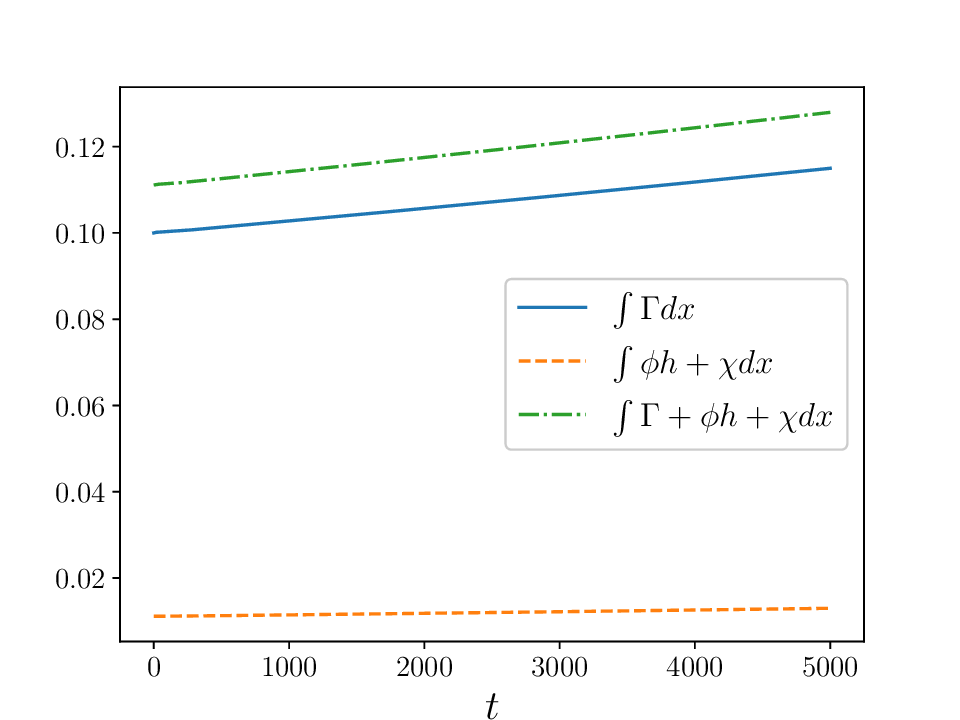}}
    \hfill
    \subcaptionbox{Mass evolution, corrected model}{\includegraphics[width=0.49\textwidth]{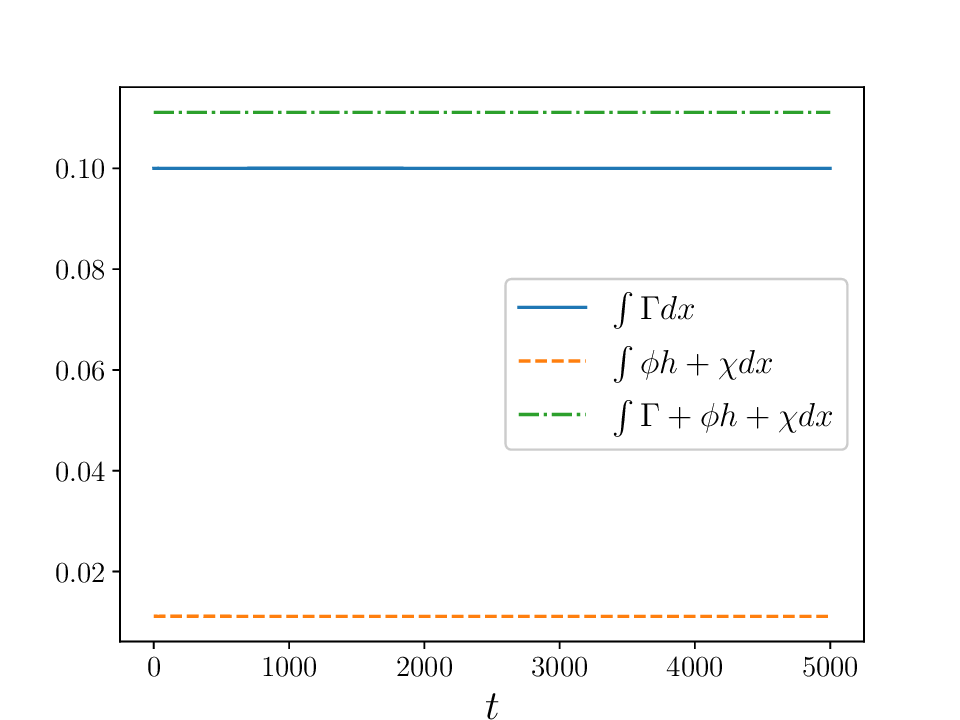}}

    \caption{
Detailed diagnostics of bulk, interfacial, and total surfactant transport for the earlier ($5/4$, left column) and corrected ($1/4$, right column) formulations. The top row shows the bulk flux contributions, the middle row shows the interfacial flux balance, and the bottom row shows the evolution of the bulk, interfacial, and total surfactant mass.
}
    \label{fig:flux_compare}
\end{figure*}

\begin{figure*}
    \centering
    \subcaptionbox{Earlier model}{\includegraphics[width=0.49\textwidth]{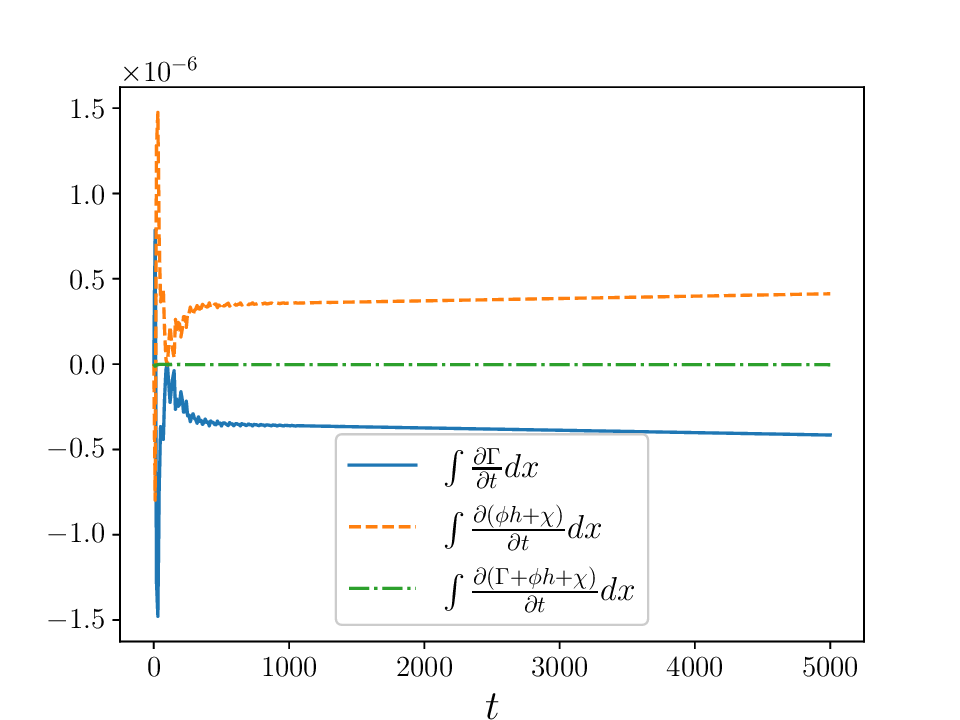}}
    \hfill
    \subcaptionbox{Corrected model}{\includegraphics[width=0.49\textwidth]{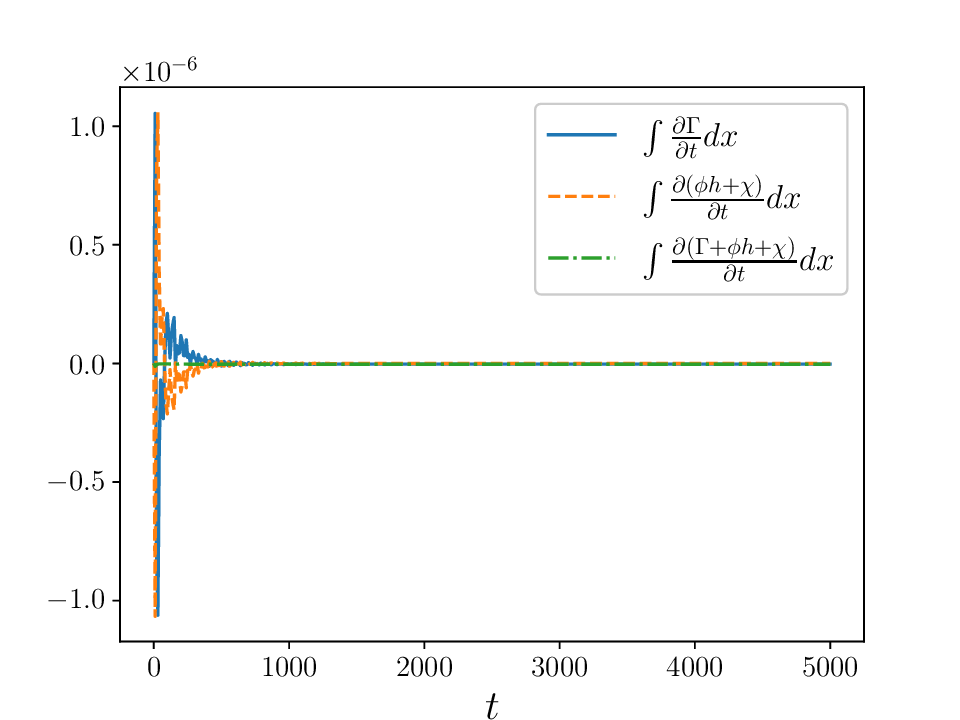}}
    \caption{
Time evolution of the rate of change of the total surfactant mass for the earlier ($5/4$) and corrected ($1/4$) formulations. The contributions from the interfacial term, bulk term, and their sum are shown.}
    \label{fig:mass_derivative_compare}
\end{figure*}
\section*{S2. Surface transport balance and its long-wave reduction} \label{sec:surface_balance}

 The nonlinear coefficient in the surface equation is fixed by the conservative structure of the surfactant balance on a moving interface. Because this point is the origin of the correction, we state the argument in a form close to Stone's derivation and then take the thin-film limit.

Let $S(t)$ be a material element of the free surface and let $\Gamma$ denote the surfactant mass per unit interfacial area. In the absence of surface diffusion and bulk exchange, Stone (1990) \cite{stoneSimpleDerivationTimedependent1990} starts from the material-surface balance
\begin{equation}
 \frac{d}{dt}\int_{S(t)} \Gamma\,dS=0 .
\label{eq:stone_integral}
\end{equation}
Using the rate of change of a material surface element, this gives \cite{stoneSimpleDerivationTimedependent1990}
\begin{equation}
 \Gamma_t+\nabla_s\cdot(\Gamma\vec{V})=0,
\label{eq:stone_full_velocity}
\end{equation}
where $\vec{V}$ is the full velocity of the moving interface and $\nabla_s=(\mathbf I-\mathbf n\mathbf n)\cdot\nabla$ is the surface-gradient operator. The surface advective flux can be decomposed into components normal to the interface, and components along the interface \cite{manikantanSurfactantDynamicsHidden2020}
$$\Gamma\vec{V} = \mathbf{n}(\mathbf{n}\cdot \vec{V}) \Gamma +(I-\mathbf{nn})\cdot \vec{V} \Gamma $$
and therefore has surface divergence \cite{manikantanSurfactantDynamicsHidden2020} 
$$\nabla_s\cdot(\Gamma\vec{V}) = \Gamma ({\boldsymbol{\nabla}_s}\cdot \mathbf{n})(\vec{V} \cdot \mathbf{n}) +\boldsymbol{\nabla}_s\cdot(\Gamma \vec{V}_s) $$
then Eq.~\eqref{eq:stone_full_velocity} becomes \cite{stoneSimpleDerivationTimedependent1990}
\begin{equation}
 \Gamma_t+ \boldsymbol{\nabla}_s\cdot(\Gamma \vec{V}_s) + \Gamma ({\boldsymbol{\nabla}_s}\cdot \mathbf{n})(\vec{V} \cdot \mathbf{n}) = 0.
\label{eq:stone_split}
\end{equation}
The last term is the surface-area dilation contribution. It is the product of the normal velocity and the mean curvature, up to the sign convention for $\mathbf n$. With surface diffusion and soluble exchange, the balance used here is therefore
\begin{equation}
 \Gamma_t+ \boldsymbol{\nabla}_s\cdot(\Gamma \vec{V}_s) + \Gamma ({\boldsymbol{\nabla}_s}\cdot \mathbf{n})(\vec{V} \cdot \mathbf{n}) = D_s \nabla_s^2 \Gamma + J_{bs}.
\label{eq:stone_soluble_split} 
\end{equation}
or equivalently Eq.~\eqref{eq:stone_full_velocity} with the full velocity inside the surface divergence and with the right-hand side added. The exchange term is not the geometric Stone term. It is the Langmuir adsorption-desorption flux
\begin{equation}
 J_{bs}=k_s\left[\kappa(1-\Gamma)C_s-\Gamma\right],
 \qquad C_s=C(x,h(x,t),t).
\label{eq:jbs_surface_section}
\end{equation}
Thus adsorption into the interface is a gain for the surface and a loss for the bulk.

We now estimate the Stone \cite{stoneSimpleDerivationTimedependent1990} area-dilation term for the thin-film graph $z=h(x,t)$. In the nondimensional variables of (\ref{dimension}), the physical slope of the interface is $\varepsilon h_x$, and the outward normal may be written as
\begin{equation}
 (\mathbf n)^*=\frac{(-\varepsilon h_x,1)}{\sqrt{1+\varepsilon^2h_x^2}} .
\label{eq:normal_graph}
\end{equation}
Restoring the dimensional scales only for the order estimate, the curvature and normal velocity are
\begin{equation}
 (\nabla_s\cdot\mathbf n)^*
 =-\frac{\varepsilon^2}{\hat h}
 \frac{h_{xx}}{(1+\varepsilon^2h_x^2)^{3/2}},
\label{eq:curvature_estimate}
\end{equation}
\begin{equation}
(\vec{V}\cdot\mathbf n)^*
 =\varepsilon\hat V\,
 \frac{w-u h_x}{\sqrt{1+\varepsilon^2h_x^2}} .
\label{eq:normal_velocity_estimate}
\end{equation}
Consequently, after nondimensionalization by $\hat h/\hat V$, Stone's normal-motion/curvature term contributes
\begin{equation}
 \frac{\hat h}{\hat V}\,\Gamma ({\boldsymbol{\nabla}_s}\cdot \mathbf{n})^*(\vec{V} \cdot \mathbf{n})^*
 =-\varepsilon^3\Gamma
 \frac{h_{xx}(w-u h_x)}{(1+\varepsilon^2h_x^2)^2} .
\label{eq:stone_term_scaled}
\end{equation}
Using the kinematic condition $w=h_t+u h_x$, this may also be written as
\begin{equation}
 -\varepsilon^3\Gamma
 \frac{h_t h_{xx}}{(1+\varepsilon^2h_x^2)^2} .
\label{eq:stone_term_kinematic}
\end{equation}
The sign depends on the direction chosen for the normal, but the order does not. The leading unsteady and advective terms in the nondimensional surface balance are $O(\varepsilon)$. Hence the explicit Stone area-dilation term is two orders smaller, namely $O(\varepsilon^3)$ in the undivided surface equation, or $O(\varepsilon^2)$ after division by $\varepsilon$. The nonlinear Marangoni correction compared in this paper appears at $O(\varepsilon Re\,Mr)$ after division by $\varepsilon$. For the standard long-wave ordering $Re=O(1)$ and $Mr=O(1)$, the area-dilation term is therefore beyond the retained order and cannot change the coefficient of the $\Gamma\Gamma_xh_x$ term.

The scalar tangential speed along the graph contains $u+\varepsilon^2wh_x$ divided by the metric factor. The $\varepsilon^2wh_x$ part again contributes only at the higher geometric order just identified. Thus, at the retained order, we use
\begin{equation}
 u_s=\frac{u(x,h(x,t),t)}{\sqrt{1+\varepsilon^2h_x^2}} .
\label{eq:us_graph}
\end{equation}
Dropping only the $O(\varepsilon^3)$ geometric correction in Eq.~\eqref{eq:stone_term_scaled}, while keeping the surface diffusion metric and the soluble exchange, gives the long-wave surface balance
\begin{widetext}
\begin{equation}
\varepsilon\left[
\Gamma_t+\frac{1}{\sqrt{1+\varepsilon^2 h_x^2}}
\frac{\partial(\Gamma u_s)}{\partial x}
\right]
=
\frac{\varepsilon^2}{Pe_s\sqrt{1+\varepsilon^2 h_x^2}}
\frac{\partial}{\partial x}
\left[
\frac{1}{\sqrt{1+\varepsilon^2 h_x^2}}\Gamma_x
\right]
+k_s\left[\kappa(1-\Gamma)C_s-\Gamma\right].
\label{eq:surface_balance_metric_v5}
\end{equation}
\end{widetext}
\section*{S3- Weighted Residual model}
We use the \citet{ruyer-quilImprovedModelingFlows2000} weighted residual technique to construct a reduced-second-order model that is valid up to a moderate Reynolds number, consistent in first order inertia, and includes second-order streamwise viscous terms.
Since the normal velocity component of the plane $w = - \int_ 0^z \partial_x u dz$ is relatively smaller than the stream component $u$, we can neglect the inertia terms and the viscous terms of the stream in the $z$ component of the momentum equation (\ref{y-momentum}). The inertia and the viscous term in the streamwise $z$ component of the momentum equation are of higher order and can therefore be omitted\cite{ruyer-quilImprovedModelingFlows2000}. The pressure distribution is obtained by integrating the remaining equation up to the order $\varepsilon$. Substituting the value of $p$ into equation (\ref{x-momentum}) produces an approximate momentum equation.
\begin{widetext}
    \begin{eqnarray} \label{reduced_mom}
    \varepsilon Re \left(\frac{\partial u}{\partial t} + u \frac{\partial u}{\partial x} + w \frac{\partial u}{\partial z}\right) = \frac{Re}{Fr^2} + \frac{\partial^2 u}{\partial z^2} 
   - \varepsilon \frac{Re}{Fr^2}\cot\theta h_x + 2 \varepsilon^2 \frac{\partial^2 u}{\partial x^2} + \varepsilon^2 \frac{\partial}{\partial x}\left[\frac{\partial u}{\partial x}\biggm|_{z = h}\right]
   + \varepsilon Ka h_{xxx},
\end{eqnarray}
\end{widetext}
where $Ka =  {\varepsilon}^2 Re We$. Let us decompose the velocity components as: $u= u_0 + \varepsilon u_1+ \varepsilon^2 u_2 + ...,\quad w = \varepsilon w_1+ \varepsilon^2 w_2 + ...$. Then equation (\ref{reduced_mom}) reduces to
\begin{widetext}
    \begin{eqnarray} \label{perturbed_mom}
    \frac{\partial^2 (u_0 + \varepsilon u_1)}{\partial z^2} + \frac{Re}{Fr^2} = \varepsilon Re \left(\frac{\partial u_0}{\partial t} + u_0 \frac{\partial u_0}{\partial x} + \varepsilon w_1 \frac{\partial u_0}{\partial z}\right) + \mathcal{K} + \varepsilon \frac{Re}{Fr^2}\cot\theta \frac{\partial h}{\partial x} 
    - \varepsilon Ka \frac{\partial^2 h}{\partial x^2} - 2 \varepsilon^2 \frac{\partial^2 u_0}{\partial x^2} 
    \nonumber\\
    - \varepsilon^2 \left[\frac{\partial u_0}{\partial x}\biggm|_h \right]  + O(\varepsilon^3) \quad
\end{eqnarray}
\end{widetext}
where $\mathcal{K}$ represents the second order inertia corrections contributed by the deviations $u_1$ from the parabolic velocity profile and $\varepsilon w_1 = - \int_0^z \partial_x u_0 dz$. Moreover, to single out $u_1$ we use the gauge condition $\int_0^h u_1 dz = 0$, which ensures that $q$ still corresponds to the local flow rate.\\
The main objective of the weighted residual method is to eliminate the dependence on $z$ by making certain assumptions about the velocity and solute concentration. We have considered the following specific profile\cite{dalessioMarangoniInstabilitiesAssociated2020} for $u$ and $C$: 
\begin{eqnarray*}
     u = \frac{Re}{Fr^2} \left(\frac{q}{h}\right) \Bar{z} \left(1 - \frac{\Bar{z}}{2}\right) + \varepsilon Re Mr \left(\frac{h}{4}\right) (2\Bar{z} - 3 \Bar{z}^2)\\
      C = \phi(x,t) + Pe_b k_s  \left(\frac{h}{2}\right)[\kappa (1 - \Gamma)\phi - \Gamma] (1 - \Bar{z}^2)
 \end{eqnarray*}
where $\Bar{z} = z/h$\\
Equation (\ref{perturbed_mom}) explains the balance of gravity acceleration and drag forces, which is influenced by inertia, viscous stream-wise diffusion, surface tension and the hydrostatic pressure gradient. Then integrating equations (\ref{mass}), (\ref{perturbed_mom}) and (\ref{concentration}) with respect to $z$ from $0$ to $h(x, t)$ and using the value of $u_s$ (the assumed velocity profile on the free surface) and $\chi$ (the excess/deficit amount of surfactant in the bulk) in the surface surfactant concentration balance equation (\ref{boundary_surface_sarfactant}) to obtain a second-order weighted residual model for the momentum equations. We chose $\Bar{z}-\Bar{z}^2/2$ as the weight function to integrate (\ref{perturbed_mom}) while unity was chosen as the weight function for the equation of the bulk concentration (\ref{concentration}). To accurately preserve the critical condition, the momentum equation must be consistent up to order $\varepsilon$ in the averaging process, as the onset of the instability is captured by the gradient expansion up to first order \citep{ruyer-quilImprovedModelingFlows2000}. After some algebra we obtain a second-order four equation weighted residual model (4-7) as given in the main text of the letter.

\section*{S4. Total surfactant mass and diagnostics}

For a soluble surfactant in a periodic domain of length $L$, the total mass at any instant $t$ is

$$M(t) = \int_0^L \int_0^h C(x,z,t) dz dx + \int_0^L  \Gamma(x,t)  dx $$
$$= \int_0^L [\phi(x,t) h (x,t) + \chi(x,t) + \Gamma(x,t) ] dx $$ 
which must remain constant for a closed system.

where, $\phi(x,t)  = \begin{cases}
C(x, z = h, t) & \text{if }\quad t > 0 \\
 \phi_e = C_e & \text{if} \quad t = 0
\end{cases}$\par
The excess/deficit amount of surfactant in the bulk is given by $$\chi = \begin{cases} \int_0^h (C - \phi)dz = \dfrac{Pe_b k_s}{3} h^2 [\kappa (1 - \Gamma)\phi - \Gamma] & \text{if}\quad t > 0 \\ (C_e- \phi_e)h_e= 0 & \text{if}\quad t=0, \end{cases}$$ where $Pe_b$ is the Peclet number, $C_e$ is the bulk concentration and $h_e$ is film thickness at the equilibrium state. The quantity $\chi$ is related to the net surfactant adsorption rate. \par
Since, the total mass of the surfactant at any instant $t$ is constant 
$$\partial_t M(t) = 0$$

\begin{eqnarray*}
     \implies 
\partial_t
\int_{0}^{L}\Bigl[\phi(x,t)\,h(x,t)+\chi(x,t)+\Gamma(x,t)\Bigr]\;dx =0,
\end{eqnarray*}

it follows that
\[
\int_{0}^{L}[
 \partial_t (\chi + \phi\,h)
\;+\;
 \partial_t \Gamma 
]\,dx =0. \]

Since \(L\) is arbitrary, the integrand must vanish:
\[
 \partial_t(\chi + \phi\,h)
\;+\;
\partial_t \Gamma =0.
\]
This is the conservation law, which provides a useful diagnostic for evaluating the consistency and fidelity of numerical models that involve soluble surfactants, which is violated by the boundary condition proposed by  Pascal et al.\cite{pascalStabilityInclinedFlow2019} and D'Alessio et al.\cite{dalessioMarangoniInstabilitiesAssociated2020}

Figures~\ref{fig:flux_compare} and \ref{fig:mass_derivative_compare} show the bulk, interfacial, and total mass-flux diagnostics for the earlier and corrected formulations. All simulations are performed for $Re = 1.5$, $Fr = 0.7071$, $Pe_b = Pe_s = 700$, $\varepsilon = 0.1$, $Mr = 1$, $k_s = 1$, $\kappa = 10$, $\Gamma_e = 0.1$, and $Ka = 0.75$ in a periodic domain of length $L = 20$. The earlier formulation exhibits a persistent imbalance between bulk and interfacial contributions, resulting in a non-zero net rate of change and a drift in the total surfactant mass. In contrast, the corrected formulation restores the balance between these contributions, with cancellations occurring to numerical accuracy, thereby ensuring global mass conservation.

\end{document}